\documentclass[a4paper,11pt]{article}

\usepackage{amsmath,amsthm,amssymb}

\renewenvironment{proof}{\noindent{\emph{\textbf{P\hspace{-0.2mm}roo\hspace{-0.15mm}f.}}}}{\hfill$\blacksquare$}

\newtheorem{theorem}{Theorem}[section]
\newtheorem{definition}[theorem]{Definition}
\newtheorem{example}[theorem]{Example}
\newtheorem{lemma}[theorem]{Lemma}
\newtheorem{proposition}[theorem]{Proposition}
\newtheorem{corollary}[theorem]{Corollary}
\newtheorem{remark}[theorem]{Remark}

\newcommand{\bs}{\boldsymbol}
\newcommand{\sr}{\scriptscriptstyle}
\newcommand{\lang}{\left\langle}
\newcommand{\rang}{\right\rangle}
\newcommand{\parfrac}[2]{\frac{\partial{#1}}{\partial{#2}}} 
\newcommand{\ad}{{\rm ad\, }}
\newcommand{\Ad}{{\rm Ad\, }}
\newcommand{\ce}{\overset{\longrightarrow}{\rm exp}\int} 
\newcommand{\e}[2]{e^{#1\,#2}}   
\newcommand{\Diff}{{\rm Diff\,}} 
\renewcommand{\Vec}{{\rm Vec\,}} 
\newcommand{\refeq}[1]{(\ref{#1})}
\newcommand{\tg}{\tilde{g\,}\!}
\newcommand{\R}{\mathbb{R}}
\newcommand{\C}{\mathcal{C}}
\newcommand{\eps}{\epsilon}
\newcommand{\la}{\lambda}
\newcommand{\Om}{\Omega}
\newcommand{\om}{\omega}
\newcommand{\up}{\Upsilon}

\newcommand{\vh}{ {\bs{\vec{h}}} }
\newcommand{\f}{\bs{f}}
\newcommand{\X}{ {\bs{X}} }
\newcommand{\Y}{\bs{Y}}
\newcommand{\Z}{ {\sr \mathcal{Z}} }
\newcommand{\coZ}{ {\sr \rm{co}\mathcal{Z}} }
\newcommand{\Plat}{{Trivializable~}}
\newcommand{\plat}{{trivializable~}}
\newcommand{\platness}{{trivializability~}}
\renewcommand{\H}{\mathcal{H}}

\newenvironment{proofLvk+bk+L2hb=0}{\noindent{\emph{\textbf{P\hspace{-0.2mm}roo\hspace{-0.15mm}f o\hspace{-0.15mm}f  P\hspace{-0.2mm}ropo\hspace{-0.2mm}s\hspace{-0.15mm}ition \ref{Lvk+bk+L2hb=0}.}}}}{\hfill$\blacksquare$}

\numberwithin{equation}{section}

\begin{document}

\author{Ulysse Serres
        \footnote{Institut \'Elie Cartan de Nancy UMR 7502, Nancy-Universit\'e/CNRS/INRIA,
                  BP 239 F-54506 Vand{\oe}uvre-l\`es-Nancy Cedex,
                  France;
                  email: \texttt{ulysse.serres@iecn.u-nancy.fr}}
}

\date{}

 \title{Control systems of zero curvature are not necessarily \plat}
 \maketitle

\begin{abstract}
A control system \(\dot{q} = \f(q,u)\) is said to be \plat if there exists local coordinates in which the system is feedback equivalent to a control system of the form
\(\dot{q} = \f(u)\).
In this paper we characterize \plat control systems and control systems for which, up to a feedback transformation,
\(\f\) and \(\partial\f/\partial{u}\) commute.
Characterizations are given in terms of feedback invariants of the system (its control curvature and its centro-affine curvature) and thus are completely intrinsic.
To conclude we apply the obtained results to Zermelo-like problems on Riemannian manifolds.
\\
\medskip

\noindent{\bf Keywords:} {Control systems, control curvature, state-feedback equivalence, Zermelo-like problems.}
\\
\medskip

\noindent{\bf MSC2000:}
34K35; 37C10; 37E35; 53B99; 93C10; 93C15
\end{abstract}

%
%
\section{Introduction}
In the present paper smooth objects are supposed to be of class $\C^\infty$.

Dynamics of the classical Riemannian geodesic problem on surfaces locally read
\begin{equation*}
\dot{q} = \cos u \,\bs{e}_1(q) + \sin u \,\bs{e}_2(q),
\quad
q \in M,
\quad
u \in S^1,
\end{equation*}
where $(\bs{e}_1,\bs{e}_2)$ is a local orthonormal frame for the Riemannian structure on $M$.
It is well known that such a system is \plat if and only if the Gaussian curvature of the surface vanishes identically.
In the Riemannian case \platness also means that all geodesics can be simultaneously rectified.

Let $M$ and $U$ be two smooth manifolds of respective dimension two and one.
Our goal in this paper is to find similar conditions for fully nonlinear two-dimensional control systems with scalar input.
More precisely, we consider smooth control systems of the type
\begin{equation}\label{contsyst}
\dot{q} = \f(q,u), \quad q \in M, \quad u \in U,
\end{equation}
whose curves of admissible velocities (or indicatrices at $q$) $u \mapsto \f(q,u)$ are strongly convex (or concave).
In other words, we suppose that \refeq{contsyst} satisfies the regularity assumptions
\begin{equation}\label{RegAss}
\bs{f}(q,u) \wedge \parfrac{\bs{f}(q,u)}{u} \neq 0,
\quad
\parfrac{\bs{f}(q,u)}{u} \wedge \parfrac{^2\bs{f}(q,u)}{u^2} \neq 0,
\quad
q \in M,\ u\in U.
\end{equation}
We aim to determine under which conditions system (\ref{contsyst}) is {\it \plat} in the sense of the following
\begin{definition}\label{def_flatness}
A control system $\dot{q}=\bs{f}(q,u)$ is said to be \plat if there exists local coordinates in which the system is (state-)feedback equivalent to a control system of
the form $\dot{q}=\bs{f}(u)$.
\end{definition}
It is worth mentioning that as well as the class of two-dimensional Riemannian manifolds, the class of two-dimensional Finsler manifolds is a particular case of control systems we treat here
(see the book \cite{ChernShen2005} as a basic and recent reference on Finsler geometry).
Indeed, if we suppose that the curves of admissible velocities of our control system are closed simple curves (in addition with the strong convexity hypothesis (\ref{RegAss})), there exists a canonically defined Finsler structure on the base manifold $M$ whose geodesics are the minimum time solutions of the optimal control problem driven by equation \refeq{contsyst}.

The Finsler analogue of the control curvature is the notion of Riemann or flag curvature (see \cite{ChernShen2005}).
Some work has been done in order to give some geometric characterizations of Finsler spaces with vanishing Riemann curvature.
In particular Mo has shown \cite{mo2001} that Finsler manifolds having zero Riemann curvature are characterized by the fact that the horizontal distribution of the projective sphere bundle has a flat foliation.
The literature on curvature problems in Finsler geometry is vast and we do not plane to observe it here.
Let us only mention the very recent paper \cite{NajafiShenTayebi2008} on the classification of Finsler metrics with scalar flag curvature
and the paper by Bao and Chern \cite{BaoChern1993} where the authors addressed the following question:
{\it describe the Finsler Spaces for which the Riemannian part of the curvature vanishes}.

Although the first part of our main result Theorem \ref{classif1} asserts that the horizontal distribution is integrable, the characterization of \plat control systems we propose here has not been treated in the framework of Finsler geometry.
For this reason we do think that the present paper can be of interest not only to people from the control theory community but also to people from the Finsler geometry community.

Section \ref{Preliminaries} contains the basics on the curvature of two-dimensional smooth control systems.
The main results are stated and proved in Section \ref{FlatSystems} which ends with some examples.
Section \ref{Conclusion} ends the paper with a summary of differences between \platness in Riemannian geometry and \platness of control systems in the sense of Definition \ref{def_flatness}
(notice this definition of \platness coincide with the definition of flatness given by Dazord in \cite{Dazord1971} for Finsler manifolds).

%
%
\section{Preliminaries}\label{Preliminaries}
In this section we briefly describe the principal feedback-invariants of two-dimensional control systems.
For more details on the subject we refer the reader to \cite{AAAbook, ulysseThese}.

\subsection{Counting the principal invariants}
Systems of the form \refeq{contsyst} are considered up to state-feedback equivalence, i.e.,
up to transformations of the form
\begin{equation}\label{eq_feedback}
\Theta(q,u)=(\phi(q),\psi(q,u)),
\end{equation}
where $\phi$ is a diffeomorphism of $M$ which plays the role of a change of coordinates 
and $\psi$ is a reparametrization of the set $U$ of controls in a way depending on the
state variable $q \in M$.
First of all, let us roughly estimate the number of parameters (invariants) in this equivalence problem.
In this case, if the coordinates on the manifold are fixed, a (germ of) control system of type (\ref{contsyst}) is parametrized
by two functions of three variables, and the group of state-feedback transformations of type (\ref{eq_feedback}) is
parametrized by two functions of two variables and one function of three variables. 
Indeed, in any local coordinate chart $q=(q_1,q_2)$ on the base manifold $M$, the control system reads
\begin{equation}\label{contsyst_coord}
\begin{split}
\dot{q}_1 & = f_1(q,u) \\
\dot{q}_2 & = f_2(q,u),
\end{split}
\end{equation}
and an element of the group of state-feedback transformations takes the form
\begin{equation*}
\Theta(q,u)=(\phi_1(q,u),\phi_2(q,u),\psi(q,u)),
\end{equation*}
where $f_1$, $f_2$, $\phi_1$, $\phi_2$ and $\psi$ are real valued functions.
Therefore, we can a priori
normalize only one function among the two functions defining control system (\ref{contsyst_coord}). Thus, we expect to have
only $2-1=1$ ``principal" feedback invariant, i.e., a function of three variables
and a certain number of feedback-invariant functions of less than three variables, in this equivalence problem.

\subsection{Curvature of two-dimensional smooth control systems}

In this section, we briefly recall some basic facts concerning the curvature of smooth control systems in dimension two.
We begin with a lemma that exhibits the duality between the set $\Vec M$ of smooth vector fields on $M$, and the set $\Lambda^1(M)$ of smooth one-forms on $M$.

We denote by $[\X,\bs{Y}]$ the Lie bracket (or commutator)
$\X \circ \bs{Y} - \bs{Y} \circ \X$ of vector fields $\X$, $\bs{Y} \in \Vec{M}$.
It is again a vector field and in local coordinates on $M$ the Lie bracket reads
$[\X,\bs{Y}](q) = (\partial{\bs{Y}}/\partial{q})\X(q) - (\partial{\bs{X}}/\partial{q})\bs{Y}(q)$.
\begin{lemma}\label{duality}
Let
$(\omega_1,\dots,\omega_n) \subset \Lambda^1(M)$ 
and 
$(\bs{f}_{1},\dots,\bs{f}_{n}) \subset \Vec M$
be two local dual basis.
Then,
$d\om_k = \sum_{i<j}c_{ij}^{k}\,\om_i \wedge \om_j$,
for all
$k \in \{1,\dots,n\}$,
if and only if
$\left[ \bs{f}_i,\bs{f}_j \right] = -\sum_{k=1}^{n}c_{ij}^{k}\bs{f}_k$,
for all $i$, $j$ $\in$ $\{1,\dots,n\}$.
\end{lemma}
%

Let us fix some notations.
For a two-dimensional smooth manifold $M$,
$\pi: T^*M \to M$ is the projection of the cotangent bundle to $M$.
We denote by $s$ the canonical Liouville one-form on $T^*M$,
$s_\la=\la \circ \pi_*$, $\la \in T^*M$.
If $\X$ is a smooth vector field on a manifold, we denote by $L_\X$ the Lie derivative along $\X$.

Consider the following time-optimal smooth control problem 
\begin{eqnarray}
&& \dot{q} = \bs{f}(q,u),\quad q \in M,\quad u \in U, \label{optpb_dyn}\nonumber \\
&& q(0) = q_0, \quad q(t_1) = q_1,     \label{optpb_bc} \\
&& t_1 \to \min \quad ({\rm or~}\max), \label{optpb_cost} \nonumber
\end{eqnarray}
where $M$ and $U$ are connected smooth manifolds of respective dimension two and one.
Denote by 
$h = \max_{u \in U}\lang \la,\bs{f}(q,u) \rang$, $\la \in T^*_qM$, $q \in M$,
the (normal) Hamiltonian function resulting from the Pontryagin Maximum Principle (PMP for short),
by $\H$ the level set $h^{-1}(\epsilon) \subset T^*M$, $\epsilon = \pm 1$,
and by $\vh$ the Hamiltonian field associated with the restriction of $h$ to $\H$.
Under the regularity assumptions of strong convexity on the curves of admissible velocities
\begin{equation}\label{Reg_Assumptions}
\bs{f}(q,u) \wedge \parfrac{\bs{f}(q,u)}{u} \neq 0,\quad
\parfrac{\bs{f}(q,u)}{u} \wedge \parfrac{^2\bs{f}(q,u)}{u^2} \neq 0, \quad
q \in M, \quad u \in U,
\end{equation}
the curve $\H_q = \H \cap T^*_qM$ admits, up to sign and translation, a natural parameter providing us with a vector field
$\bs{v}_q$ on $\H_q$ and by consequence with a vertical vector field $\bs{v}$ on $\H$.
The vector field $\bs{v}$ is characterized by the fact that it is, up to sign, the unique vector field on
$\H$ such that
\begin{equation}\label{w''=-ew+bw'}
L^2_{\bs{v}}\om = -\eps\om + b L_{\bs{v}}\om, \quad \om = s|_\H
\end{equation}
where
$b$ is a smooth
function on the level $\H$.
The function $b$, which is by definition a feedback-invariant, is called the centro-affine curvature.

The vector fields $\vh$ and $\bs{v}$ which are, by definition, feedback-invariant satisfy the nontrivial commutator relation
\begin{equation}\label{[h,[v,h]]=kv}
\Big[ \vh , \Big[ \bs{v} , \vh \Big] \Big] = \kappa \bs{v},
\end{equation}
where the coefficient $\kappa$ is defined to be {\it the control curvature} or simply {\it the curvature} of the optimal control problem (\ref{optpb_dyn})-(\ref{optpb_cost}).
The control curvature is by definition  a feedback-invariant of the control system and a function on $\H$ (and not on $M$ as the Gaussian one). Moreover, $\kappa$ is the Gaussian curvature (lifted on $\H$) if the control system defines a Riemannian geodesic problem.

From now on,
if $\theta$ denotes a parameter in the fiber such that
$\partial/\partial\theta = \bs{v}$,
we denote the Lie derivative
$L_{\bs{v}} = L_{\partial/\partial\theta} =\ '\ $.
In $(q,\theta)$ coordinates, the curvature has the form
\begin{equation}\label{curvature_form}
L_{\vh'}c - L_\vh c',
\end{equation}
where $c(q,\theta)$ is defined by
\begin{equation}\label{struct_w_w'1}
d\om = c\,\eps\om \wedge \om'.
\end{equation}
\begin{remark}{\rm
Notice that the coefficient $c(q,\theta)$ is not feedback-invariant.
Indeed, $\theta$ is only defined up to translations; in other words, any parameter
$\vartheta = \theta + \varphi(q)$ is such that
$\partial/\partial\vartheta = \partial/\partial\theta = \bs{v}$.
}
\end{remark}
This remark will be of peculiar importance in the proof of Theorem \ref{classif1}.

$(\eps\om,\om') \in \Lambda(M)$ is a coframe of differential forms on $M$ parametrized by $\theta$.
Its dual frame is $(\f,\f') \in \Vec{M}$.
The structural equations for $(\eps\om,\om',d\theta)$ are equation \refeq{struct_w_w'1} and
\begin{equation}\label{eq_struc2}
d\om' = (c'+bc)\epsilon\om \wedge \om'.
\end{equation}
The duality of the frames and Proposition \ref{duality} imply
\begin{equation}\label{[f,f']=-cf+...}
[\bs{f},\bs{f}']=-\eps c\bs{f} - (c'+bc)\bs{f}'.
\end{equation}
%

\subsection{Relation between the principal invariants of the equivalence problem}

The two relations \refeq{w''=-ew+bw'} and \refeq{[h,[v,h]]=kv} define two feedback invariants, the centro-affine curvature $b$ and the control curvature $\kappa$.
Both $b$ and $\kappa$ are functions on the three-dimensional level surface $\mathcal{H}$, so
that they are principal feedback invariants of our control system. Since our feedback equivalence problem admits
only one invariant these functions are not ``independent". Indeed we have the following proposition.
\begin{proposition}\label{Lvk+bk+L2hb=0}
The feedback invariants $b$ and $\kappa$ satisfy the following equation
\begin{equation}\label{bnk}
L_{\bs{v}}\kappa+b\kappa+L^2_{\vh}b=0.
\end{equation}
\end{proposition}
Before proving the proposition we need an auxiliary lemma.
\begin{lemma}\label{Lhb=f(c)}
Let $\theta$ be a parameter on the fiber $\H_q$ such that
$\bs{v}=\partial/\partial\theta$.
Then, the structure constant $c(\theta,q)$ defined by
$d\om = c\,\epsilon\om \wedge \om'$
satisfies
\begin{equation}\label{PDE=Lhb}
c''+bc'+\epsilon c=L_{\vh}b.
\end{equation}
\end{lemma}
\begin{proof}
Differentiating the structure equation \refeq{eq_struc2} with respect to $\theta$ leads, on the one hand, to
\begin{eqnarray*}
\left( d\om' \right)'
&=& (c''+b'c+bc')\epsilon\om \wedge \om' + (c'+bc)\epsilon\om \wedge \om'' \\
&=& (c''+b'c+bc')\epsilon\om \wedge \om' + (c'+bc)\epsilon\om \wedge (-\epsilon\om + b\,\om') \\
&=& (c''+2bc'+b'c+b^2c)\epsilon\om \wedge \om',
\end{eqnarray*}
and, taking into account \refeq{w''=-ew+bw'}, on the other hand, to
\begin{eqnarray*}
(d\om')' &=& d\om'' = d(-\epsilon\om + b\,\om')
          =  -c\,\om \wedge \om' + d_qb\,\om' + b(c'+bc)\epsilon\om \wedge \om' \\
         &=& (-\epsilon c+L_{\bs{f}}b+bc'+b^2c)\epsilon\om \wedge \om'.
\end{eqnarray*}
Summing up, we get
\begin{equation*}
c''+2bc'+b'c+b^2c=-\epsilon c+L_{\bs{f}}b+bc'+b^2c,
\end{equation*}
or equivalently,
\begin{equation*}
c''+bc'+\epsilon c=L_{\bs{f}}b-b'c=L_{\vh}b,
\end{equation*}
which ends the proof of the lemma.
\end{proof}
\\

\noindent We are now ready to prove Proposition \ref{bnk}.\\

\begin{proofLvk+bk+L2hb=0}
From lemma \ref{Lhb=f(c)} it follows that
\begin{eqnarray*}
\bigg[\parfrac{}{\theta},\bigg[\vh,\parfrac{}{\theta}\bigg]\bigg]&=&-\vh''
=-\bs{f}''+c''\parfrac{}{\theta}
=\epsilon\bs{f}+b\bs{f}'+(L_{\vh}b-bc'-\epsilon c)\parfrac{}{\theta} \\
&=&\epsilon\vh+b\vh'+L_{\vh}b\parfrac{}{\theta}.
\end{eqnarray*}
If we now compute the Lie bracket of the previous relation with $\vh$, we get for the right hand side
\begin{equation*}
\bigg[\vh,\epsilon\vh+b\vh'+L_{\vh}b\parfrac{}{\theta}\bigg]
=L_{\vh}b\vh'+b\Big[\vh,\vh'\Big]-L_{\vh}b\vh'+L_{\vh}^2b\parfrac{}{\theta}
=(b\kappa+L_{\vh}^2b)\parfrac{}{\theta},
\end{equation*}
and using Jacobi's identity, we get for the left hand side
\begin{eqnarray*}
\bigg[\vh,\bigg[\parfrac{}{\theta},\bigg[\vh,\parfrac{}{\theta}\bigg]\bigg]\bigg]
&=&-\bigg[\parfrac{}{\theta},\bigg[\bigg[\vh,\parfrac{}{\theta}\bigg],\vh\bigg]
-\bigg[\bigg[\vh,\parfrac{}{\theta}\bigg],\bigg[\vh,\parfrac{}{\theta}\bigg]\bigg] \\
&=& -\bigg[\parfrac{}{\theta},\kappa\parfrac{}{\theta}\bigg]
=\kappa'\parfrac{}{\theta}
\end{eqnarray*}
and the equation follows.
\end{proofLvk+bk+L2hb=0}\\

Notice that equation (\ref{bnk}) shows that in the special case of Riemannian problems, the curvature
$\kappa$ is a function on the base manifold $M$ without any computation. Indeed, since Riemannian problems
are characterized by the vanishing of function $b$, (\ref{bnk}) reduces to $L_{\bs{v}}\kappa=0$.

%
%
\section{\Plat systems}\label{FlatSystems}
In Riemannian geometry it is well known that if the Gaussian curvature of the surface is nonzero then, one cannot rectify simultaneously the geodesics by a change of coordinates. Only Riemannian \plat systems, i.e., systems for which the geodesics are ``straight lines" have this property.
For control systems the situation is quite different.
It is obvious that the extremals of a \plat control system
are simultaneously rectifiable
and,
that the latter implies that the curvature of the system vanishes identically,
but, in general, the converse implications do not hold.

We present here two new theorems,
the first one characterizes control systems whose extremals can be simultaneously rectified
and,
the second one characterizes \plat control systems.
Because these two characterizations are given in terms of the feedback invariants $\kappa$ and $b$,
they are intrinsic.
%

\subsection{The results}

We are now ready to state our main result:
 
\begin{theorem}\label{classif1}
There exists a feedback transformation such that
the vector fields $\bs{f}$ and $\partial\bs{f}/\partial{u}$ commute
if and only if the feedback invariants $\kappa$ and $L_{\vh}b$ are identically equal to zero.
Moreover,
if $u$ is such a parameter,
then
the infinitesimal generator of the diffeomorphism $P_u \in \Diff(M,\R^2)$ such that
\begin{equation}\label{Frob}
P_{u*}\left( \bs{f}(\cdot,u),\parfrac{\bs{f}(\cdot,u)}{u} \right)
= \left( \parfrac{}{q_1} , \parfrac{}{q_2} \right),
\end{equation}
is the vector field
\begin{equation}\label{Xu}
\bs{X}_u=(a_1(u) \pm q_2)\parfrac{}{q_1}+(a_2(u,q_2)-q_1)\parfrac{}{q_2},
\end{equation}
where the $\pm$ sign in the expression of $\X_u$ depends on whether the curves of admissible velocities of system (\ref{contsyst}) are strongly convex or strongly concave.
\end{theorem}
\begin{proof}
In this proof, we freely use the chronological calculus notation for which we refer to \cite[Chapter 2]{AAAbook}.
Suppose that $\kappa$ and $L_{\vh}b$ are identically equal to zero for control system (\ref{contsyst}).
Then, equation (\ref{[h,[v,h]]=kv}) reduces to
\begin{equation}\label{relation_courbure_nulle}
\Big[\vh,\Big[\bs{v},\vh\Big]\Big]=0.
\end{equation}
In particular, the flows $\e{t}{\vh}$ and $\e{t}{\left[ \bs{v},\vh \right]}$ commute. Therefore, the vector fields $\vh$ and $[\bs{v},\vh]$ are good candidates in order to define a system of local coordinates.
Let $\theta$ be a parameter in the fiber $\H_q$ such that
$\bs{v}=\partial/\partial\theta$.
This choice of parameter $\theta$ defines a foliation of the three-dimensional manifold $\H$, the leaves of which are formed by the trajectories of the fields $\vh$ and $\vh'$, i.e.,
\begin{equation*}
\H = \bigcup_{\la\in\H_q}\mathcal{C}_{\la},\quad
\mathcal{C}_{\la}=\left\{\e{s}{\vh'}\circ\e{t}{\vh}(\la)\ \Big|\ (s,t)\in\R^2\right\}.
\end{equation*}
Recall that this choice of $\theta$ is not feedback invariant. Indeed, the parameter $\theta$ is only fixed up to feedback transformations of the form
\begin{equation}\label{theta_plus_g} 
\theta\mapsto \pm\theta+g(q).
\end{equation}
Now fix this parameter $\theta$ in such a way that its value on the leaf $\mathcal{C}_{\la_0}$ is constant. In other words we choose the function $g$ in (\ref{theta_plus_g}) such that
\begin{equation}\label{theta|Clambda0}
\theta|_{\mathcal{C}_{\la_0}}=\theta_0.
\end{equation}
Recall that in coordinates $(\theta,q)$ on $\H$ vector fields $\vh$ and $\vh'$ take the form
\begin{equation*}
\vh = \bs{f}-c\parfrac{}{\theta}, \quad
\vh'= \bs{f}'-c'\parfrac{}{\theta},
\end{equation*}
which, in addition to (\ref{theta|Clambda0}), implies that
$c|_{\mathcal{C}_{\la_0}}$ and $c'|_{\mathcal{C}_{\la_0}}$ are zero.
Because $L_{\vh}b=0$ identically, it follows from Proposition \ref{bnk} that $c$ is solution to the Cauchy problem:
\begin{equation*}
c''+bc'+\epsilon c = 0,\quad c|_{\mathcal{C}_{\la_0}}=0,\quad c'|_{\mathcal{C}_{\la_0}}=0,
\end{equation*}
from which it follows that
$c=0$
identically on $\H$. Hence, $\vh=\bs{f}$ and $\vh'=\bs{f}'$ which, according to (\ref{relation_courbure_nulle}), is equivalent to
$[\bs{f},\partial\bs{f}/\partial{u}]=0$.
The first implication is thus proved.

We now prove the converse. Let $u$ be a control parameter such that $\f$ and
$\partial{\f}/\partial{u}$
commute. In particular,
\begin{equation*}
\bigg\langle \om,\left[\bs{f}(\cdot,u),\parfrac{\bs{f}(\cdot,u)}{u} \right]\bigg\rangle=0,
\end{equation*}
where, as usual $\om$ denotes the Liouville one-form in restriction to $\H$. According to (\ref{[f,f']=-cf+...}), one infers that
\begin{eqnarray}
0
&=& \bigg\langle \om , \left[ \bs{f} , \frac{d\theta}{du}\bs{f}' \right] \bigg\rangle
 =  \bigg\langle \om , \frac{d\theta}{du} \left[ \bs{f} , \bs{f}' \right]
                +\left( L_{\bs{f}}\frac{d\theta}{du} \right)\bs{f}'
    \bigg\rangle \nonumber \\
&=& \bigg\langle \om , \frac{d\theta}{du} \left( -\eps c\bs{f} - (c'+bc)\bs{f}' \right)
                +\left( L_{\bs{f}}\frac{d\theta}{du} \right)\bs{f}'
    \bigg\rangle \nonumber \\
&=& -c\frac{d\theta}{du} \label{0=-cdtheta/du}.
\end{eqnarray}
Because
$d\theta/du$
never vanishes (see \cite[Chapter 23, page 355]{AAAbook}), the above equation implies that $c=0$ identically on $\H$.
In this case, equations (\ref{curvature_form}) and (\ref{PDE=Lhb}) obviously imply that $\kappa$ and $L_{\vh}b$ are zero identically. The first part of the theorem is thus proved.

In order to parametrize control systems with zero curvature such that $\f$ and
$\partial{\f}/\partial{u}$
commute, we will use the classical Moser's homotopy method.
If a control system is such that
$[\f,\partial{\f}/\partial{u}]$
holds, it follows from Frobenius theorem that the vector fields $\f$ and
$\partial{\f}/\partial{u}$
can be rectified simultaneously. Thus for every $u\in U$ there exists a diffeomorphism $P_u \in \Diff M$ such that
\begin{equation}\label{Frobbis}
P_{u*}\left(\f(\cdot,u),\parfrac{\f(\cdot,u)}{u} \right) = \bigg(\f(\cdot,u_0),\parfrac{\f}{u}(\cdot,u_0)\bigg).
\end{equation}
In order to get the expression (\ref{Xu})
we use Moser's homotopy method the key idea of which is to determine the diffeomorphisms $P_u$ by representing them as the flow of a family of vector fields $\X_u$ on $M$. We thus suppose that
\begin{equation*}
\frac{d}{dt}P_u = P_u  \circ \X_u , \quad P_{u_0} ={\rm Id},
\end{equation*}
or equivalently that
\begin{equation*}
P_u= \ce_{u_0}^{u}\X_v\,dv.
\end{equation*}
The expression of $\X_u$ in coordinates will follow from the differentiation with respect to $u$ of
(\ref{Frobbis}). But, after multiplication of both sides by $P_{u*}^{-1}$, (\ref{Frobbis}) is equivalent to
\begin{eqnarray*}
\f(\cdot,u) &=& \Ad\ce_{u_0}^{u}\X_v\,dv \, \f(\cdot,u_0) , \\
\parfrac{\f(\cdot,u)}{u} &=& \Ad \ce_{u_0}^{u}\X_v\,dv \, \parfrac{\f}{u}(\cdot,u_0),
\end{eqnarray*}
which, after differentiation with respect to $u$ gives
\begin{eqnarray*}
\parfrac{\f}{u}(\cdot,u) 
&=& \Ad \ce_{u_0}^{u}\X_v\,dv \, \ad\X_u \left( \f(\cdot,u_0) \right)
  = P_{u*}^{-1} [\X_u,\f(\cdot,u_0)],  \\
\parfrac{^2\f}{u^2}(\cdot,u) 
&=& \Ad \ce_{u_0}^{u}\X_v\,dv \, \ad\X_u \left( \parfrac{\f}{u}(\cdot,u_0) \right)
  =   P_{u*}^{-1}\bigg[\X_u,\parfrac{\f}{u}(\cdot,u_0)\bigg],
\end{eqnarray*}
which, according to (\ref{Frobbis}) is equivalent to
\begin{eqnarray}
\parfrac{\f}{u}(\cdot,u_0) &=& [\X_u,\f(\cdot,u_0)], \label{Duf=adXuf} \\
P_{u*}\parfrac{^2\f}{u^2}(\cdot,u) &=& \bigg[\X_u,\parfrac{\f}{u}(\cdot,u_0)\bigg] \label{D2uf=adXuDf}.
\end{eqnarray}
Fix a system of local coordinates $q=(q_1,q_2)$ on the base manifold such that
\begin{equation*}
\f(\cdot,u_0)=\parfrac{}{q_1},\quad \parfrac{\f}{u}(\cdot,u_0)=\parfrac{}{q_2},
\end{equation*}
and denote
\begin{equation*}
\X_u=X_1(q,u)\parfrac{}{q_1}+X_2(q,u)\parfrac{}{q_2}.
\end{equation*}
In these coordinates, equation (\ref{Duf=adXuf}) reads
\begin{equation*}
-\parfrac{X_1}{q_1}=0,\quad -\parfrac{X_2}{q_1}=1,
\end{equation*}
which implies that
\begin{equation*}
X_1(q,u)=\alpha_1(q_2,u),\quad X_2(q,u)=\alpha_2(q_2,u)-q_1,
\end{equation*}
where $\alpha_1$, ans $\alpha_2$ are $C^\infty$ functions. Recall that $\f(\cdot,u)$ satisfies the second order ODE
\begin{equation*}
\parfrac{^2\f}{u^2}(\cdot,u)=-\epsilon\f(\cdot,u)-b(\cdot,u)\parfrac{\f}{u}(\cdot,u).
\end{equation*}
Thus, according to (\ref{Frobbis}), equation (\ref{D2uf=adXuDf}) reads
\begin{eqnarray*}
\bigg[\X_u,\parfrac{\f}{u}(\cdot,u_0)\bigg]
&=&  P_{u*}\parfrac{^2\f}{u^2}(\cdot,u)
 =  -\epsilon P_{u*}\f(\cdot,u) - P_{u*}\left(b(\cdot,u)\parfrac{\f}{u}(\cdot,u)\right) \\
&=& -\epsilon \f(\cdot,u_0) - b(P_u(\cdot),u) P_{u*}\parfrac{\f}{u}(\cdot,u) \\
&=& -\epsilon \f(\cdot,u_0) - b(P_u(\cdot),u) \parfrac{\f}{u}(\cdot,u_0).
\end{eqnarray*}
So in our system of local coordinates on $M$ this last equation reads
\begin{equation*}
- \parfrac{X_1}{q_2} = - \parfrac{\alpha_1}{q_2}= - \epsilon , \quad
- \parfrac{X_2}{q_2} = - \parfrac{\alpha_2}{q_2}= - b(P_u(q),u),
\end{equation*}
from which it follows that
\begin{equation*}
X_1(q_1,q_2)=a_1(u)+\epsilon q_2,\quad X_2(q_1,q_2)=a_2(q_2,u)-q_1,
\end{equation*}
which is the required expression for the field $\X_u$ and ends the proof.
\end{proof}\\
\begin{remark}{\rm
Notice that commutativity between vector fields $\bs{f}$ and
$\partial\bs{f}/\partial{u}$
is not a feedback-invariant property.
When the curvature is identically zero the Theorem \ref{classif1} shows that equation $L_\vh b=0$ reduces to the nonautononous ODE
$dq/du=\bs{X}_u(q)$.
}
\end{remark}
\begin{remark}{\rm
Using the variation formula described in \cite[Chapter 2, Section 2.7]{AAAbook}, one easily sees that the diffeomorphism $P_u$ takes the form
\begin{equation*}
P_u = \ce_{u_0}^{u}
      e^{(s-u_0){{\phantom{-}0\ \pm 1}\choose{-1\ \phantom{\pm}0}} } \textstyle{a_1(s)\phantom{q_2,} \choose a_2(q_2,s)}ds
      \circ e^{(u-u_0) {{\phantom{-}0\ \pm 1}\choose{-1\ \phantom{\pm}0}} }.
\end{equation*}
}
\end{remark}

The following theorem characterizes \plat control systems.
\begin{theorem}\label{classif2}
A control system of type (\ref{contsyst}) is \plat if and only if its feedback invariants $\kappa$, $L_{\vh}b$ and $L_{[\boldsymbol{v},\vh]}b$ vanish identically.
\end{theorem}
\begin{proof}
Suppose that the system under consideration is \plat. By definition this system is feedback equivalent to a system of the form $\dot{q}=\bs{f}(u)$. For such a system it is obvious that the feedback invariant $b$ depends only on the control parameter $u$ and that the Hamiltonian is horizontal. Therefore, the feedback invariants $\kappa$, $L_{\vh}b$ and  $L_{[\bs{v},\vh]}b$ vanish identically.

We now prove the converse. It follows from Theorem \ref{classif1} that the vanishing of $\kappa$ and $L_{\vh}b$ implies that, up to a feedback, the vector fields ${\vh}$ and  $[\bs{v},\vh]$ are horizontal. Therefore, the vanishing of $L_{\vh}b$ and  $L_{[\bs{v},\vh]}b$ is equivalent to the vanishing of $L_{\bs{f}}b$ and  $L_{[\bs{v},\bs{f}]}b$, from which it immediately follows that the invariant $b$ depends only on the control parameter $u$. In this case, the infinitesimal generator of the one-parameter family of diffeomorphisms defined by (\ref{Frob}) is
\begin{equation*}
\bs{X}_u=(a_1(u) \pm q_2)\parfrac{}{q_1}+(a_2(u)-q_1)\parfrac{}{q_2}.
\end{equation*}
Thus,
\begin{equation*}
\left( \ce_{u_0}^{u}\bs{X}_v\,dv \right)(q)
= \e{(u-u_0)}{{{\phantom{-}0\ \pm 1}\choose{-1\ \phantom{\pm}0}}}(q)
+\int_{u_0}^{u}\e{(u-v)}{{{\phantom{-}0\ \pm 1}\choose{-1\ \phantom{\pm}0}}} \textstyle{a_1(v) \choose a_2(v)}dv
\end{equation*}
from which it follows that
\begin{eqnarray*}
\bs{f}(q,u) = \Ad  \ce_{u_0}^{u}\bs{X}_v\,dv  \, {1 \choose 0}= \bs{f}(u).
\end{eqnarray*}
That ends the proof.
\end{proof}\\
%

\subsection{Examples}
In this section $(M,g)$ denotes a two-dimensional Riemannian manifold and $(\bs{e}_1,\bs{e}_2)$ denotes a local $g$-orthonormal frame.
\begin{example}
{\rm \textbf{\Plat Riemannian manifolds}.\\
Both Theorems \ref{classif1} and \ref{classif2} imply the following classical theorem.
\begin{theorem}
A two-dimensional Riemannian manifold is \plat if and only if its Gaussian curvature vanishes identically.
\end{theorem}
In the Riemannian case, the control curvature is the Gaussian curvature (see e.g. \cite{AAAbook, ulysseThese}).
Moreover, in this case the feedback invariant $b$ vanishes identically, which shows that Theorems \ref{classif1} and \ref{classif2} lead to the same thesis. If we denote by $(\bs{e}_1,\bs{e}_2)$ a local orthonormal basis for the Riemannian structure on the manifold, we then see that Theorems \ref{classif1} and \ref{classif2} reduce to
\begin{equation*}
\kappa\equiv 0\quad
\Leftrightarrow
\quad\textrm{there exists a feedback such that }\ \bigg[\bs{f},\parfrac{\bs{f}}{u}\bigg]=[\bs{e}_1,\bs{e}_2]=0.
\end{equation*}
On the other hand, if $[\bs{e}_1,\bs{e}_2] = 0$, according to the Frobenius theorem, one can find a system of local coordinates on $M$ such that
$\bs{e}_1=\partial/\partial{q_1}$, $\bs{e}_2=\partial/\partial{q_2}$,
i.e., such that the dynamics of the Riemannian problem read
\begin{equation*}
\dot{q}=\cos u\parfrac{}{q_1}+\sin u\parfrac{}{q_2}=\bs{f}(u).
\end{equation*}
Consequently, the system is \plat.
Moreover, if we fix local coordinates on the base manifold and set $b=0$ in the proof of Theorem \ref{classif1} we see that \plat Riemannian problems are parametrized by the vector field
\begin{equation*}
\bs{X}_u=q_2\parfrac{}{q_1}-q_1\parfrac{}{q_2}.
\end{equation*}
}
\end{example}
%

\begin{example}\label{ChapClassifZerExamp}
{\rm
\textbf{\Plat Zermelo-like problems.}
Zermelo's navigation problem on a two-dimensional Riemannian manifold $(M,g)$ is the time-optimal control problem:
\begin{eqnarray*}
&& \dot{q}=\bs{X}(q)+u,\quad q \in M,\quad u \in S^1, \\
&& q(0)=q_0,\quad q(t_1)=q_1                          \\
&& t_1\to \min,                                       \\
\end{eqnarray*}
The restriction to $h^{-1}(1)$ of the Hamiltonian vector field associated to the maximized Hamiltonian function $h$
reads (see \cite[Chapter 3]{ulysseThese})
\begin{equation*}
\vh(q,u) = \X(q) + \cos u\,\bs{e}_1 + \sin u\,\bs{e}_2 - c_\Z(q,u)\parfrac{}{u},
\end{equation*}
where
\begin{eqnarray}\label{cZ_formula}
c_\Z(q,u)
 &=& \cos^2 u\,L_{\bs{e}_2}X_1 + \cos u\sin u (L_{\bs{e}_2}X_2 - L_{\bs{e}_1}X_1) + \sin^2 u\,L_{\bs{e}_1}X_2  \nonumber \\
 & & + (1 + \cos u\,X_1 + \sin u\,X_2)(c_1\cos u + c_2\sin u), \\
 & & \nonumber \\
X_1(q) &=& \lang \X(q),\bs{e}_1(q) \rang_g, \quad X_2(q) = \lang \X(q),\bs{e}_2(q) \rang_g. \nonumber
\end{eqnarray}
Suppose that a Zermelo's navigation problem Theorem \ref{classif1}.
Hence, according to the proof of the same theorem, there exists a vertical parameter $\vartheta$ such that
$\partial/\partial\vartheta = \bs{v}$
and $[\bs{f},\bs{f}'] = 0$.
For any smooth functions $f$, $g$ and any smooth vector fields $\X$, $\Y$, the general relation
$[f\X , g\Y] = fg[\X ,\Y] + f L_{\X}g\Y - g L_{\Y}f\X$
and an easy calculation imply that vector fields $\f$ and
$\partial\f/\partial{u}$
satisfy the nontrivial commutation relation
\begin{equation*}
\Big[ \f , \parfrac{\f}{u} \Big] = -\frac{\eps c_\Z}{\varphi} \f + \alpha \parfrac{\f}{u} +\beta\parfrac{}{u},
\quad \alpha, \ \beta \in \C^\infty(\H).
\end{equation*}
Then, a similar computation as the one made to obtain \refeq{0=-cdtheta/du} shows that there exists a system of local coordinates on $M$
such that $c_\Z(q,u)$ equals to zero identically.
In particular, in such a coordinates system we have
\begin{equation*}
2c_1 = c_\Z(q,0)-c_\Z(q,\pi) = 0, \quad 2c_2 = c_\Z(q,\pi/2)-c_\Z(q,-\pi/2) = 0,
\end{equation*}
which shows that the Riemannian manifold must be \plat. If we choose local coordinates on $M$ in which $\bs{e}_1$, $\bs{e}_2$ commute,
then, according to (\ref{cZ_formula}), the vanishing of $c_\Z$ implies in particular that
\begin{eqnarray*}
 2L_{\bs{e}_1} X_1 &=& c_\Z(q,\pi/4) + c_\Z(q,-\pi/4) = 0, \quad L_{\bs{e}_2} X_1 = c_\Z(q,0) = 0, \\
 2L_{\bs{e}_1} X_2 &=& c_\Z(q,\pi/4) - c_\Z(q,-\pi/4) = 0, \quad L_{\bs{e}_2} X_2 = c_\Z(q,\pi/2) = 0,
\end{eqnarray*}
which trivially implies that the coordinates $X_1$, $X_2$ of the drift in $(\bs{e}_1,\bs{e}_2)$ have to be constant.

Summing up, we have proved the following
\begin{theorem}
A Zermelo navigation problem on a Riemannian manifold is \plat if and only if the Riemannian manifold is \plat and the drift vector field is constant in any system of local coordinates in which $\bs{e}_1$ and $\bs{e}_2$ commute.
\end{theorem}
We now turn our attention to the co-Zermelo problem for which we refer the reader to \cite{ulysseco-Zer} for details.
Let $\up$ be a one-form on $M$ such that $|\up|_g < 1$.
We call co-Zermelo problem of the pair $(g,\up)$ the following time-optimal control problem on $M$
\begin{eqnarray*}
&& \dot{q}={\displaystyle\frac{u}{1 + \lang \up_q , u \rang}},\quad q\in M,\quad u\in T_qM,\quad |u|_g=1,\\
&& q(0)=q_0,\quad q(t_1)=q_1, \\
&& t_1\to \min.
\end{eqnarray*}
The Hamiltonian function of PMP is of this problem reads
\begin{equation*}
h(\la)
= \frac{ - \lang \la , \up_{\pi(\la)} \rang_g 
         + \sqrt{\lang \la , \up_{\pi(\la)} \rang_g^2 + \left( 1 - |\up_{\pi(\la)}|_g^2 \right) |\la|_g^2}}
{1 - |\up_{\pi(\la)}|_g^2}.
\end{equation*}
Let $(\bs{e}_1,\bs{e}_2)$ be a $g$-orthonormal frame and parametrize fibers by $u$ in such a manner that
$\lang \la - \up_{\pi(\la)},\bs{e}_1(\pi(\la)) \rang = \cos\theta$,
$\lang \la - \up_{\pi(\la)},\bs{e}_2(\pi(\la)) \rang = \sin\theta$.
Then, the Hamiltonian field takes the form
\begin{equation*}
\vh(q,u)
= \frac{1}{\varphi(q,u)} \left( \cos u\,\bs{e}_1(q) + \sin u\,\bs{e}_2(q) 
                                    +(c_g(q,u) + \Om(q))\parfrac{}{u}
                              \right),
\end{equation*}
where $\Om \in \C^\infty(M)$ and $\varphi \in \C^\infty(\H)$ are the functions defined by
\begin{equation*}
d\up = -\Om\,dV_g, \quad \varphi(q,u) = 1 + \cos u\lang \up_q,\bs{e}_1(q) \rang + \sin u\lang \up_q,\bs{e}_2(q) \rang,
\end{equation*}
and the curvature reads
\begin{equation}\label{Kco-Zer}
\kappa_\coZ^{\sr (g,\up)}
= \varphi^{-2}\left( \kappa_g + \Om^2 + \sin uL_{\bs{e}_1}\Om - \cos uL_{\bs{e}_2}\Om
  -\mathcal{S}(\varphi) \right),
\end{equation}
where $\mathcal{S}(\varphi)$, the Schwartzian derivative of $\varphi$ is defined by
\begin{equation*}
\mathcal{S}(\varphi) = \varphi L_{\vh}\left( \frac{L_{\vh}\varphi}{2} \right)
                      -\left( \frac{L_{\vh}\varphi}{2} \right)^2.
\end{equation*}

An easy computation and \cite[Proposition 3.4 and Corollary 3.5]{ulysseco-Zer},
(which assert that a given co-Zermelo problem on $(M,g)$ is feedback-equivalent to a Zermelo problem on the same manifold $M$ equipped with a Riemannian metric $\tg$ generally different from $g$) imply the following
\begin{corollary}
A co-Zermelo problem on a Riemannian manifold is \plat if and only if the Riemannian manifold is \plat and the drift one-form is constant in any system of local coordinates in which $\bs{e}_1$ and $\bs{e}_2$ commute.
\end{corollary}

}\end{example}
%

%
%
\section{Conclusion}\label{Conclusion}
The differences, exhibited by Theorems \ref{classif1} and \ref{classif2}, between the Riemannian and the control cases
({\bf RC} and {\bf CC} respectively) can be summarized as follows:
\begin{displaymath}
\begin{array}{lccccc}
\textrm{{\bf RC}:}\quad&
\textrm{\plat}&~\Leftrightarrow~&~[\f,\f'] = 0 \textrm{~(up to a feedback)}~&~\Leftrightarrow~&~\kappa \equiv 0 \\
\textrm{{\bf CC}:}\quad&
\textrm{\plat}&~\Rightarrow    ~&~[\f,\f'] = 0 \textrm{~(up to a feedback)}~&~\Rightarrow    ~&~\kappa \equiv 0 
\end{array}
\end{displaymath}
We want to point out that the existence of a feedback such that extremals project onto $M$ as straight lines neither implies that the control system is \plat,
nor the existence of a feedback such that $\f$ and
$\partial\f/\partial{u}$
commute.
It can be easily seen if one considers a co-Zermelo problem on the Euclidean plane $\R^2$ whose drift is an exact form $\up = df$.
In this case the Hamiltonian field $\vh$ dynamic on $M$ reads (see \cite{ulysseco-Zer})
\begin{equation*}
\vh = \frac{1}{\varphi(q,u)} \left( \cos u\parfrac{}{q_1} + \sin u\parfrac{}{q_2} \right),
\quad
\varphi(q,u) = 1 + \cos u\parfrac{f}{q_1} + \sin u\parfrac{f}{q_2},
\end{equation*}
and, according to \refeq{Kco-Zer}, the curvature reads
$\kappa=-\mathcal{S}(\varphi)$, which has no reason to be identically zero.
Indeed, one can check that taking $f = q_1^2 + q_2^2$ leads to
$\kappa = 3(1 + 2q_1\cos u + 2q_2\sin u)^{-4}$.
%


\begin{thebibliography}{1}

\bibitem{AAAbook}
A.~A. Agrachev and Y.~L. Sachkov.
\newblock {\em Control theory from the geometric viewpoint}, volume~87 of {\em
  Encyclopaedia of Mathematical Sciences}.
\newblock Springer-Verlag, Berlin, 2004.
\newblock Control Theory and Optimization, II.

\bibitem{BaoChern1993}
D.~Bao and S.-S. Chern.
\newblock On a notable connection in {F}insler geometry.
\newblock {\em Houston J. Math.}, 19(1):135--180, 1993.

\bibitem{ChernShen2005}
S.-S. Chern and Z.~Shen.
\newblock {\em Riemann-{F}insler geometry}, volume~6 of {\em Nankai Tracts in
  Mathematics}.
\newblock World Scientific Publishing Co. Pte. Ltd., Hackensack, NJ, 2005.

\bibitem{Dazord1971}
P.~Dazord.
\newblock Tores finsl\'eriens sans points conjugu\'es.
\newblock {\em Bull. Soc. Math. France}, 99:171--192; erratum, ibid. 99 (1971),
  397, 1971.

\bibitem{mo2001}
X.~Mo.
\newblock Finsler spaces with vanishing {R}iemann curvature.
\newblock {\em Rev. Roumaine Math. Pures Appl.}, 46(4):455--463, 2001.

\bibitem{NajafiShenTayebi2008}
B.~Najafi, Z.~Shen, and A.~Tayebi.
\newblock Finsler metrics of scalar flag curvature with special
  non-{R}iemannian curvature properties.
\newblock {\em Geom. Dedicata}, 131:87--97, 2008.

\bibitem{ulysseThese}
U.~Serres.
\newblock {\em {G}\'eom\'etrie et classification par feedback des syst\`emes de
  contr\^ole non lin\'eaires de basse dimension}.
\newblock PhD thesis, Universit\'e de Bourgogne, Dijon, $24{\rm th}$ March
  2006.
\newblock In English. Available at:\\
  {\texttt{http$:$//tel.archives-ouvertes.fr/tel-00172902/en/}}.

\bibitem{ulysseco-Zer}
U.~Serres.
\newblock On {Z}ermelo-like problems: {G}auss-{B}onnet inequality and {E}.
  {H}opf theorem.
\newblock {\em J. Dynam. Control Systems}, 15(1):99--131, 2009.

\end{thebibliography}


\end{document}